\documentclass[11pt,a4paper,reqno]{amsart}
\usepackage{amsaddr, amsmath, amsfonts,amsthm, amssymb, bm, hyperref, booktabs, algorithm, algorithmic, subcaption, graphicx}

\usepackage[backend=biber, style=numeric, giveninits=true, isbn=false, doi=true
]{biblatex}
\DefineBibliographyStrings{english}{in = {}}
\usepackage{siunitx}
\usepackage[a4paper,margin=2cm]{geometry}

\title[KL expansion of TRBB and Inference of Unit Distribution]{Karhunen-Lo\`{e}ve Expansion of Time-Rescaled Brownian Bridge and Inference of Unit Distribution}
\date{}
\graphicspath{ {Figures/} }
\newcommand{\diff}{\mathop{}\!\mathrm{d}}
\let\mathscr\relax
\DeclareMathAlphabet{\mathscr}{U}{rsfs}{m}{n}
\numberwithin{equation}{section}
\begin{document}

\author{Leonardo Chhabra}
\address{G2Q Computing SRL, Italy}

\author{Anindya Goswami}
\address{Department of Mathematics, IISER Pune, India}

\author{Nikhil Thekkedath Haridas}
\address{Department of Mathematics, IISER Pune, India}
\email{nikhilth271828@gmail.com}

\theoremstyle{plain}
\newtheorem{theorem}{Theorem}[section]
\newtheorem{lemma}[theorem]{Lemma}
\newtheorem{proposition}[theorem]{Proposition}
\newtheorem{corollary}[theorem]{Corollary}
\theoremstyle{definition}
\newtheorem{definition}[theorem]{Definition}
\newtheorem{example}[theorem]{Example}
\theoremstyle{remark}
\newtheorem{remark}[theorem]{Remark}

\numberwithin{equation}{section} 

\begin{abstract}
We study the Karhunen–Lo\`{e}ve (KL) basis of the time-rescaled Brownian bridge (TRBB) with time rescales of the form $t\mapsto \int_0^t f^2(s)\diff s$. For testing unit distributions, we construct a new statistic based on the projection of empirical process onto the KL basis functions of the TRBB and derive its asymptotic distribution. Moreover, using this statistic, we develop a goodness-of-fit test for uniformity and study its power against a specific non-uniform alternative. Numerical experiments indicate that, in the settings considered, the proposed test can have higher power than classical alternatives such as the Kolmogorov–Smirnov and Cram\'{e}r–von Mises tests. Furthermore, we propose a two-parameter family of TRBBs and derive the closed-form KL expansion for this new deterministically time-changed process. The family includes a diverse class of monotone rate time rescales, both increasing and decreasing, while retaining a tractable eigenstructure. Finally, an analysis with real data reveals greater flexibility of the proposed parametric family for modelling monotone densities on $[0,1]$.
\end{abstract}

\maketitle

\section{Introduction}
Karhunen-Loève (KL) expansion of a stochastic process $\{X(t):t\in I\}$ is analogous to the Fourier series of functions on $I$, where $I$ is a compact interval in $\mathbb{R}$. The process is represented as the coordinate expansion of an orthonormal basis in $L^2(I)$ consisting of eigenfunctions of its covariance operator. In contrast to the Fourier series, the basis here is dependent on the process. This makes the expansion optimal, as it achieves the least possible mean-squared error for any finite truncation compared to other basis expansions (see \cite[Thm.~7.2.8]{hsing}). Moreover, the expansion coefficients are random variables rather than fixed numbers and are uncorrelated. Thus, the KL expansion is a \textit{biorthogonal} expansion. For a Gaussian process, the coefficients are normally distributed, implying that they are independent. 

The derivation of KL expansion of various stochastic processes has been carried out, but not limited to, by the following papers: \cite{Maccone1986,  timechangedt^a03, klmcBM06, klmcBB07, klgenexp08, alphabridge10, kl1detreBM12, kladditive2013,  kloubridge14, kl2detreBM14,  klmdetreBM14, KLlevy18, }. KL expansion is important for various applications including, but not limited to, a) Quantization of Gaussian processes \cite{2002, hal}, b) Large and small deviation asymptotics \cite{timechangedt^a03, alphabridge10,  kladditive2013, kl1detreBM12, kl2detreBM14, klmdetreBM14}, c) Functional data analysis \cite{chiou2004functional, hsing, hu2022sparse}, and d) Goodness-of-fit tests. In particular, KL expansion of the Brownian bridge has a natural application in statistical inference. In \cite{timechangedt^a03}, KL expansions for the broad family $t^\beta B(t^\alpha) \;\;(\alpha>0)$ are derived, from which the unweighted time-rescaled Brownian bridge (TRBB) $B(t^\alpha)$ follows as a special case. The eigenstructure is expressed through Bessel functions and its zeros. Their analysis yields exact distributional formulas and tail asymptotics for the associated $L^2$-norms, and they point out the connection of these norms with weighted goodness-of-fit statistics, including Cram\'{e}r–von Mises (CvM) type tests. A similar statistical use of KL expansions appears in \cite{klgenexp08}, which is about a Gaussian process generated by independent exponential pairs, where the explicit KL decomposition is used to analyze a CvM type test for bivariate independence. 

Comparison between the empirical distribution function (EDF) and the cumulative distribution function (CDF) leads to a class of goodness-of-fit statistics, known as EDF statistics \cite{Darling1952, DarlingonKSCvM1957, Durbin1973, Stephens1974}. These include Kolmogorov-Smirnov (KS), CvM, Anderson-Darling (AD), Kuiper, and Watson statistics. In \cite[Ch.~5]{shorack2009empirical}, the KL expansion of centered Gaussian processes determined by their covariance function is used to represent statistics of integral tests, such as the CvM and AD statistics, as infinite weighted sums of chi-square random variables. The principal components of the empirical process have been used to find an orthogonal decomposition of the CvM statistic in \cite{DurbinKnott1972, DurbinKnottTaylor1975}. The above is merely a representative collection, and a comprehensive literature survey is out of the scope of this paper.

For testing unit distributions, we construct a new statistic based on the projection of empirical process onto the principal components (KL basis functions) of TRBBs, and derive its asymptotic distribution. Moreover, we consider a goodness-of-fit test for uniformity using this statistic, derive its mean and variance under a specific non-uniform alternative, and study the corresponding asymptotic power. We then compare the power of the proposed test with that of the classical KS and CvM tests. The numerical experiments illustrate the higher power in the settings considered.

In \cite{Maccone1986}, general time-rescaled Brownian motion is studied under time rescales of the form $t\mapsto \int_0^t f^2(s)\diff s$. In that setting, the eigenfunctions are expressed in terms of Bessel functions of the first kind, with the order depending on $t$. This dependence makes their computation more expensive and can increase numerical error when evaluating test statistics. While certain specific choices of time rescaling may alleviate these issues, such restrictive choices can limit practical applicability. This motivates the search for a parametric family of time rescales that is both sufficiently rich and admits a tractable eigenstructure. In this work, we study a TRBB motivated by applications to unit distributions and monotone densities on $[0,1]$. We introduce a two-parameter family of time rescales with monotone (both increasing and decreasing) rates and derive an explicit KL expansion for the corresponding TRBB. Specifically, for $t\in[0,1]$, we consider the function
\begin{equation}\label{generalf}f(t) =(at+b-\min (0,a))^{-\frac{2}{3}}, \quad b>0,  \end{equation}
where $a=0$ gives the Brownian bridge with no variable scaling. On the other hand, positive and negative values of $a$ give decreasing and increasing rates, respectively. To the best of our knowledge, this class of time rescales has not been studied previously. Furthermore, we show that \eqref{generalf} provides a flexible and competitive model for monotone densities on $[0,1]$.

The paper is organized as follows. Section~\ref{notationsec} introduces the notation and basic preliminaries used throughout the paper. In Section~\ref{testsec}, we define the new test statistic, called projection statistic, and derive its asymptotic distribution. In Section~\ref{powersec}, we obtain the asymptotic power of the proposed test for uniformity against an alternative. In Section~\ref{KLTRBBsec}, we derive the KL expansion of the proposed two-parameter family of TRBBs. Section~\ref{numericalsec} presents the eigenfunctions, variance behavior, and power comparison. In Section~\ref{datasec}, we fit the density corresponding to \eqref{generalf} to a real dataset and compare the goodness-of-fit statistics with benchmark models. Finally, we add some concluding remarks to Section~\ref{conlsec}.

\section{Notations}\label{notationsec}
Let $(\Omega, \mathcal{F},\mathbb{P})$ be a probability space and $I$ a compact interval in $\mathbb R$. Let $X=\{X(t): t \in I\}$ be a measurable second-order process ($\mathbb{E}[X^2(s)] <\infty, \;\forall s \in I$) with sample paths almost surely in $L^2(I)$.  We assume $X$ is a centered process, meaning $\mathbb{E}[X(t)] = 0$ for all $t \in I$. The covariance function $k_X(s, t):= \mathbb E\, [X(s) X(t)]$ satisfies $ \int_I k_X(s, s) \diff s<\infty$. Then $X$ can be seen as an $L^2(I)$-valued random vector (see \cite[Prop.~3.19]{DaPrato_Zabczyk_2014}) with 
$$\mathbb E [\|X\|^2_{L^2(I)} ] =\int_Ik_X(s, s) \diff s<\infty.$$
Since for all $f \in L^2(I)$, $\left\langle\mathbb{E}[X],f\right\rangle =0$,  $X$ is a centered $L^2(I)$-valued random vector. For $g,h\in L^{2}(I)$, the covariance operator $\mathscr{C}_X: L^2(I) \to L^2(I)$ is defined by
\begin{align*} \mathscr{C}_Xg :=& \mathbb{E}[\langle X,g\rangle X], \\ \langle\mathscr{C}_Xg,h\rangle=& \mathbb{E}[\langle X,g\rangle \langle X,h\rangle].\end{align*}
Note that $\mathscr{C}_X$ is self-adjoint, non-negative definite, and compact. The eigenvalue-eigenfuction pairs are denoted by $\{\lambda_k , \phi_k\}_{k \geq 1}$, where $\lambda_{1}\geq\lambda_{2}\geq\cdots>0$ are positive eigenvalues. By the spectral theorem, the corresponding eigenfunctions $\phi_k\in L^{2}(I)$ form a complete orthonormal system (CONS) for the closure $\overline{\mathrm{Im}(\mathscr{C}_X)}$. For a second-order process, the covariance operator takes the form
$$(\mathscr{C}_X y)(\cdot)=\int_I y(s) k_X(s, \cdot) \diff s, \quad \forall y \in L^2(I).$$ 
 
By KL theorem \cite{KAR1946, kacsiegert1947, loeve1948}, if $k_X$ is continuous on $I\times I$, then
\begin{equation}\label{klconvergence} \lim_{N\rightarrow \infty}\sup_{t\in I}\; \mathbb{E} \left[ \Bigl( X(t) - \sum_{k=1}^{N} x_k \phi_k(t) \Bigr)^{\!2} \,\right]
=0, \end{equation}
where coefficients $x_k$ are given by $x_k(\omega)=\int_I X(t,\omega) \phi_k(t) \diff t$, $\forall \omega \in \Omega$, satisfying $\mathbb E[x_k]=0$ and $\mathbb E [x_k x_j]= \lambda_k\delta_{k j}$. If $\xi_k=\frac{x_k}{\sqrt{\lambda_k}}$, then 
\begin{equation} \label{klexpansion}X(t)=\sum_{k=1}^{\infty} \sqrt{\lambda_k }\xi_k \phi_k(t),\end{equation} 
where the convergence is in $L^2(\Omega)$ and $\xi_k$ are centered, mutually uncorrelated random variables with unit variance. By \eqref{klconvergence}, the convergence in \eqref{klexpansion} is uniform on $I$. If the process is Gaussian, then $\xi_i \overset{i.i.d.}{\sim} \mathcal{N}(0,1)$. The proof of the KL theorem relies on Mercer's theorem for the decomposition of real-valued kernels \cite{Mercer1909}. A generalization of the theorem to Hilbert space-valued stochastic processes is stated and proved in \cite{KLH}.

Let $(\Omega, \mathcal{F}, \{\mathcal{F}_t\}_{t\geq0},\mathbb{P})$ be a filtered probability space with the adapted Brownian motion $\{W(t):t\geq 0\}$. For $c>0$, the Brownian bridge process on $[0,c]$ is defined as 
\begin{equation} \label{BB} B^{c}(t):=W(t)-\frac{t}{c}W(c),\quad \forall t \in [0,c].  \end{equation} 
Therefore, $B:= B^1$ denotes the standard Brownian bridge on $[0,1]$. As $B^c$ is a centered process, the covariance function is given by  
\begin{equation}\label{CFBB}  k_{B^c}(s,t)=\mathbb E [B^c(s)B^c(t)]=\min (s,t)-\frac{st}{c},\end{equation}
for all $t,s\in[0,c]$ (see \cite[Ch.~5, Lem.~6.9]{karatzas}). Let $f:[0,1]\to \mathbb{R}$ be such that 
\begin{equation}\label{c}c:=\int_0^1f^2(s)\diff s < \infty.   \end{equation} 
Following \cite{Maccone1986}, we define a TRBB as 
\begin{equation} \label{TRBB}Y(t) := B^c\left(\int_0^tf^2(s)\diff s\right),\quad \forall  t \in [0,1], \end{equation}
 where $B^c$ is the Brownian bridge on $[0,c]$ defined by \eqref{BB},  and $c$ is given by \eqref{c}. For all $t_1,t_2 \in [0,1] $, the covariance function of $Y$ takes the form 
\begin{align} k_Y(t_1, t_2) := \operatorname{Cov}(Y(t_1), Y(t_2))&= \min \left(\int_0^{t_1}f^2\diff t,\int_0^{t_2}f^2\diff t\right)-\frac{1}{c}\int_0^{t_1}f^2\diff t\int_0^{t_2}f^2\diff t \label{CVfnTRBB}\\ &= \int_0^{t_1\land t_2}f^2\diff t -\frac{1}{c}\int_0^{t_1}f^2\diff t\int_0^{t_2}f^2\diff t, \nonumber \end{align} 
which is continuous on $[0,1]^2$. Consequently, the covariance operator $\mathscr{C}_Y: L^2([0,1]) \to L^2([0,1])$ is written as
\begin{equation}\label{COTRBB}\mathscr{C}_Y g(t)= \int_0^1 k_Y(s,t) g(s)\diff s, \quad \forall g\in  L^2([0,1]),\quad \forall t \in [0,1].\end{equation}
\section{Unit distribution and Goodness-of-fit}\label{testsec}
In this section, we discuss goodness-of-fit tests for absolutely continuous CDFs where the support has a bounded and connected interior. The test statistic is constructed using the eigenfunctions and eigenvalues of the covariance operator of a TRBB. For the sake of simplicity, we restrict the discussion to unit distributions \cite{Unitdistributionreview}. Therefore, $f$ is chosen so that $\int_0^1 f^2(s)\diff s =1$. Let $0\le T_1,\dots,T_n\le 1$ be a random sample from a unit distribution. We consider the null hypothesis that the CDF of that distribution is given by 
$$F(x) = \begin{cases}0, & x<0,\\ \int_0^xf^2(s)\diff s, & 0\leq x \leq 1, \\ 1, & x>1.\end{cases}$$ 
The EDF $F_n$ is defined as 
$$F_n(t):=\frac{1}{n}\sum_{i=1}^n\mathbf{1}_{\{T_i\le t\}},$$ 
where $\mathbf{1}_{\{\cdot\}}$ is the indicator function. Subsequently, the empirical process \cite{EPtextbook2023, shorack2009empirical} is defined as 
$$\alpha_n(t) := \sqrt{n}(F_n(t)-F(t)).$$ 
Here we consider the TRBB 
\begin{equation}\label{UDTRBB}\{B(F(t)):{t\in [0,1]}\}.\end{equation} 
Let $\{\varLambda_k,\varphi_k\}_{k\ge1}$ be the eigenpairs of its covariance operator $\mathscr{C}_{B(F)}$, sorted so that $\varLambda_1\geq\varLambda_2\geq\dots>0$.
\begin{theorem}\label{theo4.1} For each $k\ge1$, define the coefficient of projection 
\begin{equation}\label{eq:gamma-def}\gamma_k^{(n)}:=\int_0^1 \alpha_n(t)\varphi_k(t)\diff t.\end{equation}
Then $\gamma_k^{(n)}$ can be written as
\begin{equation}\label{eq:gamma-iid-sum}\gamma_k^{(n)}=\frac{1}{\sqrt{n}}\sum_{i=1}^n\eta_k(T_i),\qquad\text{where}\quad\eta_k(u):=\int_u^1\varphi_k(t)\diff t-\int_0^1 F(t)\varphi_k(t)\diff t.\end{equation}
Moreover, under the null hypothesis $T_i\overset{i.i.d.}{\sim} F$, for any fixed $m\ge1$, the vector $\bigl(\gamma_1^{(n)},\dots,\gamma_m^{(n)}\bigr)$ converges in distribution, as $n$ tends to $\infty$, to the Gaussian vector $\bigl(Y_1,\dots,Y_m\bigr)$ with mean zero and covariance
\begin{equation}\label{eq:cov-eta}\mathrm{Cov}(Y_j,Y_k) =\int_0^1\eta_j(u)\eta_k(u)f^2(u)\diff u=\varLambda_k\,\delta_{jk}.\end{equation}
\end{theorem}
\begin{remark}
The empirical process $\alpha_n$ converges in distribution to $B(F)$ in the Skorokhod space \cite{shorack2009empirical}. This result is known as Donsker's theorem \cite{Donsker1952} in empirical process theory, and it should imply the result stated in Thm.~\ref{theo4.1}. Nevertheless, we provide a proof here using the central limit theorem (CLT).
\end{remark}
Define the random variable 
\begin{equation}\label{TRBBstatistic}s_n :=\frac{\sum_{j=1}^m \gamma_j^{(n)}}{\sqrt{\sum_{j=1}^m\varLambda_j}}.\end{equation}
Under the null hypothesis that the sample is drawn from the distribution $F$, the test statistic $s_n$ - which we call the projection statistic - asymptotically follows a standard normal distribution. Specifically, by applying Thm.~\ref{theo4.1} and letting $n \to \infty$, we obtain
$$s_n=\frac{\frac{1}{\sqrt{n}}\sum_{i=1}^n\sum_{j=1}^m\eta_j(T_i)}{\sqrt{\sum_{j=1}^m\varLambda_j}}\overset{d}{\longrightarrow} \mathcal{N}(0,1)$$
\begin{proof}[Proof of Theorem~\ref{theo4.1}]
By definition of $\alpha_n$ we have
\begin{align*}\gamma_k^{(n)}&=\sqrt{n}\int_0^1\bigl(F_n(t)-F(t)\bigr)\varphi_k(t)\diff t.\end{align*}
Since $F_n(t)=\tfrac{1}{n}\sum_{i=1}^n\mathbf{1}_{\{T_i\le t\}}$, we may exchange sum and integral to obtain
\begin{align*}\int_0^1 F_n(t)\varphi_k(t)\diff t &= \frac{1}{n}\sum_{i=1}^n\int_0^1\mathbf{1}_{\{T_i\le t\}}\varphi_k(t)\diff t=\frac{1}{n}\sum_{i=1}^n\int_{T_i}^1\varphi_k(t)\diff t.\end{align*}
Consequently
\begin{align*}\gamma_k^{(n)} &= \sqrt{n}\Bigl(\frac{1}{n}\sum_{i=1}^n\int_{T_i}^1\varphi_k(t)\diff t - \int_0^1 F(t)\varphi_k(t)\diff t\Bigr)\\&=\frac{1}{\sqrt{n}}\sum_{i=1}^n\Bigl(\int_{T_i}^1\varphi_k(t)\diff t-\int_0^1 F(t)\varphi_k(t)\diff t\Bigr) =\frac{1}{\sqrt{n}}\sum_{i=1}^n\eta_k(T_i),\end{align*}
which proves \eqref{eq:gamma-iid-sum}. 
For a fixed $m$, the random vectors
$$\bigl(\eta_1(T_i),\dots,\eta_m(T_i)\bigr),\qquad i=1,\dots,n,$$
are i.i.d. with finite second moments since $\varphi_k \in L^2([0,1])$. Clearly, under the null hypothesis, the expectation of $\eta_k(T_i)$ is zero for all $i=1,\ldots, n$. Therefore, the multivariate CLT applies to the sum in \eqref{eq:gamma-iid-sum}, giving
\begin{equation*}\frac{1}{\sqrt{n}} \sum_{i=1}^n(\eta_1(T_i),\dots,\eta_m(T_i)) \overset{d}{\longrightarrow} \mathcal N(0,\Sigma),\end{equation*}
under the null hypothesis, where $\Sigma$ is the covariance matrix with entries
$$\Sigma_{jk}=\mathrm{Cov}\bigl(\eta_j(T_i),\eta_k(T_i)\bigr)=\int_0^1\eta_j(u)\eta_k(u)f^2(u)\diff u.$$
Thus $(\gamma_1^{(n)},\dots, \gamma_m^{(n)})$ converges to a mean-zero Gaussian vector with covariance given by the above integral. In order to calculate $\Sigma_{jk}$, we write $\mu^F_k:=\int_0^1 F(t)\varphi_k(t)\diff t$ so that $\eta_k(u)=\int_u^1\varphi_k(t)\diff t-\mu^F_k$. Since we have
\begin{align*}\int_0^1\Bigl(\int_u^1\varphi_j(t)\diff t\Bigr)\mu^F_k\;f^2(u)\diff u&= \mu^F_k\int_0^1\Bigl(\int_u^1\varphi_j(t)\diff t\Bigr)f^2(u)\diff u \\ &= \mu^F_k\int_0^1\left(\int_0^t f^2(u)\diff u\right)\varphi_j(t)\diff t\\&= \mu^F_k\mu^F_j,\end{align*}
we get,
\begin{align*}\Sigma_{jk}=\int_0^1\eta_j(u)\eta_k(u)f^2(u)\diff u &= \int_0^1\Bigl(\int_u^1\varphi_j(t)\diff t\Bigr)\Bigl(\int_u^1\varphi_k(s)\diff s\Bigr)f^2(u)\diff u - \mu^F_j\mu^F_k \\&= \int_0^1\int_0^1\varphi_j(t)\varphi_k(s)\Bigl(\int_{0}^{\min(t,s)}f^2(u) \diff u\Bigr)\diff t\diff s -\mu^F_j\mu^F_k \\&=\int_0^1 \int_0^1\varphi_j(t)\varphi_k(s)\min\left(\int_{0}^{t}f^2(u) \diff u,\int_{0}^{s}f^2(u) \diff u\right)\diff t\diff s -\mu^F_j\mu^F_k.\end{align*}
Using \eqref{CVfnTRBB}, $k^F(s,t):= k_{B(F(\cdot))} = \min(F(s),F(t)) - F(s)F(t)$. Thus
\begin{equation*}\Sigma_{jk}=\int_0^1\int_0^1\varphi_j(t)\varphi_k(s)k^F(t,s)\diff t\diff s.\end{equation*}
Using the eigenvalue equation
$$\int_0^1 k^F(t,s)\varphi_k(s)\diff s=\varLambda_k\varphi_k(t)$$
and orthonormality of the $\varphi_k$'s, we conclude
\begin{align*}\Sigma_{jk} &= \int_0^1\varphi_j(t)\Bigl(\int_0^1 k^F(t,s)\varphi_k(s)\diff s\Bigr)\diff t\\&=\varLambda_k\int_0^1\varphi_j(t)\varphi_k(t)\diff t \\ &= \varLambda_k\,\delta_{jk}.\end{align*}
This proves \eqref{eq:cov-eta} and completes the proof of the theorem.  
\end{proof}
\section{Test for uniformity}\label{powersec}
It is worth noting that the projection statistic $s_n$, introduced in the earlier section, is defined using the arbitrarily chosen unit distribution $F$ mentioned in the null hypothesis. While the earlier section derives the asymptotic distribution of $s_n$ under the null hypothesis, here we obtain the asymptotic distribution of $s_n$ under an alternative as well. Although the presentation involves uniformity as null, it has a trivial extension to an arbitrary unit null distribution. In \cite{DurbinKnottTaylor1975}, the authors discuss the asymptotic powers of three test statistics against alternatives of interest. Here, we have the null and alternative hypotheses for testing uniformity of unit distributions as follows. For all $t\in [0,1]$
\begin{align}\mathrm H_0 &: \;\; f(t)=1,  \label{H_0}\\ \label{H_A} \mathrm H_A &: \;\;  f(t)= g(t),\end{align}
where $g^2$ is a known density function. The CDF corresponding to \eqref{H_0} and \eqref{H_A} are $R(t):=t$ and $S(t):= \int_0^t g^2(s)\diff s$ respectively. We refer to the test using $s_n$ \eqref{TRBBstatistic} as the projection test. Let $\{\lambda^B_k,\phi^B_k\}_{k\ge1}$, $\lambda^B_1\geq \lambda^B_2\geq \dots>0,$ be the eigenpairs of the covariance operator of standard Brownian bridge on $[0,1]$. Let $T_1,\dots,T_n$ be a sample from $[0,1]$. In hypothesis testing, under the null hypothesis \eqref{H_0}, the expression of $\gamma_k^{(n)}$ \eqref{eq:gamma-iid-sum} is 
\begin{equation}\label{neweta}\gamma_k^{(n)}=\frac{1}{\sqrt{n}}\sum_{i=1}^n\eta_k(T_i),\qquad\text{where}\quad\eta_k(u):=\int_u^1\phi^B_k(t)\diff t-\int_0^1 R(t)\phi^B_k(t)\diff t.\end{equation}
Now, under $H_A$, we find the expectation of $\eta_k(T_i)$ to be 
\begin{equation}\label{expunderalter}\mathbb{E}_A [\eta_k(T_i)]= \int_0^1 S(t)\phi^B_k(t)\diff t-\int_0^1 R(t)\phi^B_k(t)\diff t= m_{1,k}-m_{2,k},    \end{equation} 
where $m_{1,k}:= \int_0^1 S(t)\phi^B_k(t)\diff t$ and $m_{2,k}:= \int_0^1 R(t)\phi^B_k(t)\diff t$. The covariance under the alternative is given by
$$\operatorname{Cov}_A(\eta_j(T_i), \eta_k(T_i))= \int_0^1\eta_j(u)\eta_k(u)g^2(u)\diff u-(m_{1,j}-m_{2,j})(m_{1,k}-m_{2,k}), \quad 1\leq j,k \leq m.$$
Now using \eqref{neweta}, 
\begin{align*}&\int_0^1\eta_j(u)\eta_k(u)g^2(u)\diff u \\&=  \int_0^1 \int_0^1\mathbf{1}_{[u,1]}\phi^B_j(t)\diff t \int_0^1\mathbf{1}_{[u,1]}\phi^B_k(s)\diff s \;g^2(u)\diff u  -m_{1,j}m_{2,k} - m_{1,k}m_{2,j} +m_{2,j}m_{2,k} \\&= \int_0^1\int_0^1\phi^B_j(t)\phi^B_k(s)\Bigl(\int_{0}^{\min(t,s)}g^2(u) \diff u\Bigr)\diff t\diff s -m_{1,j}m_{2,k} - m_{1,k}m_{2,j} +m_{2,j}m_{2,k} \\&=\int_0^1 \int_0^1\phi^B_j(t)\phi^B_k(s)\min\left(\int_{0}^{t}g^2(u) \diff u,\int_{0}^{s}g^2(u) \diff u\right)\diff t\diff s-m_{1,j}m_{2,k} - m_{1,k}m_{2,j} +m_{2,j}m_{2,k}. \end{align*}
Consequently, using \eqref{CVfnTRBB}, 
$$\operatorname{Cov}_A(\eta_j(T_ i), \eta_k(T_i)) = \int_0^1 \int_0^1\phi^B_j(t)\phi^B_k(s)\left(k^S(s,t)+S(t)S(s)\right)\diff t\diff s -m_{1,j}m_{1,k},$$
where $k^S$, the covariance function of $\{B(S(t)): t\in [0,1]\}$, is $k^S(s,t)=\min(S(s),S(t)) - S(s)S(t)$.
Thus, 
\begin{equation}\label{covunderalternate}\operatorname{Cov}_A(\eta_j(T_i), \eta_k(T_i)) = \int_0^1 \int_0^1\phi^B_j(t)\phi^B_k(s)k^S(s,t)\diff t\diff s = \left\langle \mathscr{C}_{B(S)}\phi_j^B, \phi_k^B\right\rangle_{L^2[0,1]},\end{equation}
where $\mathscr{C}_{B(S)}$ is the covariance operator. 

For brevity, we introduce the notations
$$M := \sum_{j=1}^m (m_{1,j}-m_{2,j}), \quad V:=  \sum^m_{j,k=1}\left\langle \mathscr{C}_{B(S)}\phi_j^B, \phi_k^B\right\rangle_{L^2[0,1]}, \quad \Delta_n:=\frac{\sqrt{n}\,M}{\sqrt{\sum_{j=1}^m\lambda_j^B}}, \quad \sigma:=\sqrt{\frac{V}{\sum_{j=1}^m\lambda_j^B}}.$$
From \eqref{neweta}, \eqref{expunderalter} and \eqref{covunderalternate}, it follows that $\sum_{j=1}^m \gamma_j^{(n)}$ has mean $\sqrt{n}\,M$ and variance $V$ under $H_A$. Furthermore, under $H_A$, using the strong law of large numbers, the projection statistic \eqref{TRBBstatistic} satisfies
$$\frac{s_n}{\sqrt{n}}=\frac{\sum_{j=1}^m \gamma_j^{(n)}}{\sqrt{n\sum_{j=1}^m\lambda^B_j}} \longrightarrow \frac{M}{\sqrt{\sum_{j=1}^m\lambda_j^B}}.$$
Now we want to estimate the power of the projection test. The power of the test is the probability of correctly rejecting the null hypothesis. A two-sided level $\alpha$ test would reject $H_0$ when $|s_n|>z_{\alpha/2}$. If $n$ is very large, then using CLT, under $H_A$, we may write 
$$\frac{\frac{s_n}{\sqrt{n}}-\frac{\Delta_n}{\sqrt{n}}}{\frac{\sigma}{\sqrt{n}}} \approx \mathcal{N}(0,1).$$
Subsequently, we may write $s_n \approx \mathcal{N}(\Delta_n, \sigma^2)$. We denote the power of projection test as $\text{power}_{\text{proj}}(n)$. Thus, the asymptotic power is calculated via the standard normal CDF $\Phi$ as
$$\text{power}_{\text{proj}}(n)=1-\Big(\Phi\Big(\frac{z_{\alpha/2}-\Delta_n}{\sigma}\Big)-\Phi\Big(\frac{-z_{\alpha/2}-\Delta_n}{\sigma}\Big)\Big),$$
which converges to $1$ as $n$ tends to $\infty$. However, it is important to estimate the power for small values of $n$. To this end, we discuss the following steps. Fix a grid of sample sizes $n$, projection dimension $m$, significance $\alpha$, simulation sizes $K_{\text{null}}$ and $K_{\text{alt}}$. The algorithm is given below.
\begin{algorithm}
\renewcommand{\thealgorithm}{}
\caption{Empirical power calculation using finite sample Monte Carlo}
\label{Poweralgo}
\begin{algorithmic}[1]
\STATE $K_{\text{null}}$ number of random samples of size $n$ are drawn under $H_0$ (i.e., $U[0,1]$).
\STATE For each $l=1,\ldots, K_{\text{null}}$, $s_n^{(l)}$, i.e., the $s_n$ value for $l$-th sample, is computed.
\STATE Form the finite-sample two-sided critical value
$$c_n \;=\; \text{quantile}_{1-\alpha}\bigl(\,\lvert s_n^{(1)}\rvert,\dots,\lvert s_n^{(K_{\text{null}})}\rvert\,\bigr),$$
the value below which $100(1-\alpha)\%$ of the absolute test statistics fall. 
\STATE $K_{\text{alt}}$ number of random samples of size $n$ are drawn under $H_A$, and $\bar s_n^{(l)}$, the $s_n$ value for $l$-th sample, is computed.
\STATE Compute the empirical power
$$\widehat{\text{power}}(n)=\frac{1}{K_{\text{alt}}}\sum_{l=1}^{K_{\text{alt}}}\mathbf{1}_{\{\,|\bar s_n^{(l)}|>c_n\,\}}.$$
\STATE Clearly, it is an estimator of the success probability of a Bernoulli distribution, and its standard error estimate is given as
$$\widehat{\mathrm{SE}} = \frac{1}{\sqrt{K_{\text{alt}}}}\sqrt{\frac{K_{\text{alt}}\;\widehat{\text{power}}(n)\bigl(1-\widehat{\text{power}}(n)\bigr)}{K_{\text{alt}}-1}} \approx \sqrt{\frac{\widehat{\text{power}}(n)\bigl(1-\widehat{\text{power}}(n)\bigr)}{K_{\text{alt}}}}.$$
\end{algorithmic}
\end{algorithm}
In the next section, we discuss the derivation of the eigenfunctions and eigenvalues of the covariance operator of a TRBB.
\section{KL expansion of TRBB} \label{KLTRBBsec} 
In addition to $f\in L^2([0,1])$, we further assume its continuity and $f((0,1))\subset (0,\infty)$ in this section. In order to derive the KL expansion of the process $Y$, given in \eqref{TRBB}, we recall the expression of the covariance operator \eqref{COTRBB}, which is in the simplified form of an integral operator with kernel \eqref{CVfnTRBB}, so that further analysis can be carried out on the resulting eigenvalue equation. The covariance operator is self-adjoint, non-negative definite, and compact. Therefore, the eigenfunctions of $\mathscr{C}_Y$, $\{\phi_{k}\}_{k\geq 1}$, form a CONS for $L^2([0,1])$ with the corresponding non-negative eigenvalues $\{\lambda_{k}\}_{k\geq 1}$. From \eqref{TRBB}, $\{Y(t): t\in [0,1]\}$ is a Gaussian process, and \eqref{klexpansion} gives
$$Y(t)= \sum_{k=1}^\infty \sqrt{\lambda_k} \xi_k \phi_k(t), \quad t \in [0,1],$$ 
where $\xi_k \overset{i.i.d.}{\sim} \mathcal{N}(0,1)$. The eigenvalue equation is 
\begin{equation} \label{eveqn} \int_0^1 \left[\int_0^{t_1\land t_2}f^2\diff t\right]\phi_k(t_2)\diff t_2 - \frac{1}{c}\int_0^1 \left[\int_0^{t_1}f^2\diff t\int_0^{t_2}f^2\diff t\right]\phi_k(t_2)\diff t_2 = \lambda_k\phi_k(t_1). \end{equation}
From now on, we will consider $\lambda_k>0$. The eigenfunctions corresponding to positive eigenvalues are uniformly continuous (see \cite[Lem.~3.2]{KLH}). The first term on the left side of \eqref{eveqn} can be rewritten as the sum of integrals on intervals $[0,t_1]$ and $[t_1,1]$, where the integrands have continuous partial derivatives in $t_1$. Hence, using Leibniz's integral rule, the left side of \eqref{eveqn} is differentiable w.r.t. $t_1$. Therefore, $\phi_k$ is differentiable for each $k$. Also, from \eqref{eveqn}, 
\begin{equation}\label{bcode} \phi_k(0) = 0, \quad \phi_k(1) = 0.\end{equation} 
So, we have two boundary conditions required to solve the upcoming ODE. Differentiating \eqref{eveqn} w.r.t. $t_1$ yields 
$$f^2(t_1)\left[ \int_{t_1}^1 \phi_k (t_2)\diff t_2 + \frac{1}{c}\int_0^1 \left[\int_0^{t_2}f^2\diff t\right]\phi_k(t_2)\diff t_2 \right]= \lambda_k \phi_k^\prime(t_1), \quad t_1 \in (0,1),$$
where $\phi_k '$ denotes the first order derivative of $\phi_k $. Using $f>0$ on $(0,1)$, we get 
$$\int_{t_1}^1 \phi_k (t_2)\diff t_2 + \frac{1}{c}\int_0^1 \left[\int_0^{t_2}f^2\diff t\right]\phi_k(t_2)\diff t_2 = \lambda_k \frac{\phi_k^\prime(t_1)}{f^2(t_1)}.$$
As before, it can again be argued that the left side of above is differentiable w.r.t. $t_1$. Consequently, differentiating both sides w.r.t $t_1$, we see that \eqref{eveqn} reduces to the following Sturm-Liouville ODE
\begin{equation}\label{SLode}\frac{\diff}{\diff t} \left[ \frac{1}{f^2(t)} \phi_k'(t) \right] + \frac{1}{\lambda_k} \phi_k(t) = 0,\end{equation}
with boundary conditions \eqref{bcode}. In order to express the solution to \eqref{SLode}, using some elementary functions, we consider some specific forms of $f$. In particular, we consider
\begin{equation}\label{f} f_1(t):=(at+b)^{-\frac{2}{3}} , \quad \text{and} \quad  \quad f_2(t):=(a(1-t)+b)^{-\frac{2}{3}},\end{equation} 
where $a,b>0$. An example is plotted in Fig.~\ref{f1f2}. It is worth noting that $\int_0^1 f_1^2(s)\diff s$ and $ \int_0^1 f_2^2(s)\diff s$ are equal. We denote this scalar as $\tilde c$, which depends on the parameters $a$ and $b$. The processes corresponding to the pair of functions in \eqref{f} are denoted as 
\begin{equation} \label{the2TRBB} Y_1(t) := B^{\tilde{c}}\left(\int_0^tf_1^2(s)\diff s\right), \quad \text{and} \quad Y_2(t) := B^{\tilde{c}}\left(\int_0^tf_2^2(s)\diff s\right), \quad t \in [0,1].   \end{equation} 
Both time rescalings are included in the following unified representation 
\begin{equation} \label{model1} t\mapsto \int_0^t (as +b - \min (0,a)\;)^{-\frac{4}{3}}\diff s = \begin{cases}  \int_0^t (as+b)^{-\frac{4}{3}}\diff s, & \; a\geq0,\\\int_0^t (|a|(1-s)+b)^{-\frac{4}{3}}\diff s, & \; a<0,\end{cases}\end{equation} 
where $b>0$. This produces the time for $Y_1$, $Y_2$, and $B^{b^{-\frac{4}{3}}}$ when $a$ is positive, negative, and zero, respectively. In the remainder of this section, we derive the KL expansion for the first two processes. Following \cite{Maccone1986}, we wish to express the eigenfunctions in terms of eigenvalues. To this end, we propose the following results.
\begin{theorem}\label{mainthm}
Let $f(t) =(at+b-\min (0,a))^{-\frac{2}{3}}$ for $a\neq0$ and $b>0$. The eigenfunctions $\phi_k$ with positive eigenvalues $\lambda_k$ satisfy \eqref{SLode} with boundary conditions \eqref{bcode}. Define the following for $t\in [0,1]$:    
$$\chi(t) = \sqrt{f(t) \int_0^t f(s)  \diff s},  \quad \psi_k(t) = \frac{1}{\sqrt{\lambda_k}} \int_0^t f(s)  \diff s.$$
Furthermore, set $y_k(t) = \frac{\phi_k}{\chi}(t)$ for $t\in (0,1]$. Then $y_k$ satisfies the ODE 
\begin{equation}\label{psiode} \psi_k(t)\frac{\diff}{\diff t} \left[ \frac{\psi_k(t)}{\psi^\prime_k(t)} y^\prime_k(t) \right] + \left[ \psi^2_k(t) - \left(\frac{1}{2}\right)^2 \right] \psi^\prime_k(t) y_k(t) = 0, \quad t\in (0,1). \end{equation}
\end{theorem}
\begin{proof}
For $t\in (0,1)$, ODE \eqref{SLode} takes the form 
$$\frac{\chi(t)}{f^2(t)}\,y^{\prime\prime}_k(t) + \left( \frac{\chi'(t)}{f^2(t)} + \frac{\diff}{\diff t}\!\Big(\frac{\chi(t)}{f^2(t)}\Big) \right) y^{\prime}_k(t) + \left( \frac{\diff}{\diff t}\!\Big(\frac{\chi'(t)}{f^2(t)}\Big) + \frac{\chi(t)}{\lambda_k} \right) y_k(t) = 0 .$$
Multiplying both sides of the above equation by $\chi(t)$ gives
\begin{equation}\label{adjoint}\frac{\diff}{\diff t}\!\Bigg( \frac{\chi^2(t)}{f^2(t)}\,y^{\prime}_k(t) \Bigg)+ \chi(t)\left( \frac{\diff}{\diff t}\!\Big(\frac{\chi'(t)}{f^2(t)}\Big) + \frac{\chi(t)}{\lambda_k} \right) y_k(t) = 0.\end{equation}
In order to show that \eqref{adjoint} is the same as \eqref{psiode}, we must have \begin{equation}  \label{****}\chi(t)\left( \frac{\diff}{\diff t}\!\Big(\frac{\chi'(t)}{f^2(t)}\Big) + \frac{\chi(t)}{\lambda_k} \right)=\left[ \psi^2_k(t) - \left(\frac{1}{2}\right)^2 \right] \frac{\psi^\prime_k(t)}{\psi_k(t)},\end{equation}
since $\frac{\psi_k(t)}{\psi^\prime_k(t)} = \frac{\chi^2(t)}{f^2(t)}$. Using the fact that 
$$\frac{\chi^2(t)}{\lambda_k}= \psi_k^2(t)\frac{\psi^\prime_k(t)}{\psi_k(t)},\;\; \forall k\geq 1 ,$$
and $\frac{\psi^\prime_k(t)}{\psi_k(t)}$ is independent of $k$, \eqref{****} is equivalent to 
\begin{equation} \label{egworkshere}\chi(t) \frac{\diff}{\diff t}\!\Big(\frac{\chi'(t)}{f^2(t)}\Big)=-\frac{\psi^\prime_k(t)}{4\psi_k(t)}, \quad \forall t \in (0,1). \end{equation}
Hence, proving \eqref{psiode} reduces to showing \eqref{egworkshere} when $f=(at+b-\min(0,a))^{-\frac{2}{3}}$. Let
$$I(t):=\int_0^t f(s)\diff s,\qquad A(t):=\chi(t)^2=f(t)I(t).$$
Note that $A' = f'I + f^2$ and $A'' = f''I + 3 f f'$. Since $\chi=\sqrt{A}$, we have
$$\chi'=\frac{A'}{2\chi},\qquad\chi''=\frac{A''}{2\chi}-\frac{(A')^2}{4\chi^3}.$$
Define
$$\mathcal{E}(t):=\chi(t)\frac{\diff}{\diff t}\!\Big(\frac{\chi'(t)}{f^2(t)}\Big)=\chi(t)\Big(\frac{\chi''}{f^2}-2\frac{\chi' f'}{f^3}\Big)(t)=\left(\frac{\chi\chi''}{f^2}-2\frac{\chi\chi' f'}{f^3}\right)(t).$$
Substituting the expressions for $\chi\chi''$ and $\chi\chi'$ gives
$$\mathcal{E}=\frac{A''}{2f^2}-\frac{A'^2}{4Af^2}-\frac{A' f'}{f^3}.$$
Now, by replacing $A=fI$, $A'=f'I+f^2$, $A''=f''I+3ff'$, the three terms are
\begin{align*}\frac{A''}{2f^2}&= \frac{f''I+3ff'}{2f^2}= \frac{f''I}{2f^2} + \frac{3f'}{2f},\\[4pt]-\frac{A'^2}{4Af^2}&= -\frac{(f'I+f^2)^2}{4f^3I}= -\frac{(f')^2I}{4f^3} -\frac{f'}{2f} -\frac{f}{4I},\\[4pt]-\frac{A' f'}{f^3}&= -(f'I+f^2)\frac{f'}{f^3}= -\frac{(f')^2 I}{f^3} -\frac{f'}{f}.\end{align*}
Summing these three contributions yields
\begin{align}\label{mathcal{E}}\mathcal{E}= \frac{f''I}{2f^2} - \frac{5}{4}\frac{f'^2 I}{f^3} - \frac{f}{4I}= \frac{I}{4f^3}\big(2f f'' - 5 f'^2\big) - \frac{f}{4I}.\end{align}
For $f$ given in \eqref{generalf} and $g(t):=at+b-\min(0,a)$, we have
$$f=g^{-2/3},\qquad f'=-\frac{2a}{3}g^{-5/3},\qquad f''=\frac{10a^2}{9}g^{-8/3},$$
and
$$2f f'' = \tfrac{20a^2}{9}g^{-10/3} = 5 f'^2,$$
hence $2f f''-5f'^2=0$. Substituting this into \eqref{mathcal{E}}, the expression reduces to $\mathcal{E} = -\frac{f}{4I}$, i.e.
$$\chi(t)\,\frac{\diff}{\diff t}\!\Big(\frac{\chi'(t)}{f(t)^2}\Big)= -\frac{f(t)}{4\displaystyle\int_0^t f(s)\diff s},$$
which is the claimed identity \eqref{egworkshere}.
\end{proof}
Now we make the observation that \eqref{psiode} takes the form of a modified Bessel equation. We make it explicit with the following corollary.
\begin{corollary}\label{cormainthm} Let $y_k$ and $\psi_k$ be as in Thm.~\ref{mainthm}. For all $k\geq 1$, we have $$y_k(t) =A_kJ_{\frac{1}{2}}(\psi_k(t)),\quad t\in [0,1], \;\;A_k\in \mathbb R,$$ where $J_\nu$ is the Bessel function of the first kind of order $\nu$.
\end{corollary}
\begin{proof}
As $\psi_k$ is monotonic, we make the substitution
\begin{equation*}x = \psi_k(t),\quad t\in [0,1],\end{equation*} 
and observe that \eqref{psiode} becomes the Bessel equation
\begin{equation*}x\frac{\diff}{dx} \left[ x \frac{d\Psi_k(x)}{dx} \right] + \left( x^2 - \left(\frac{1}{2}\right)^2 \right) \Psi_k(x) = 0,\quad x\in (0,\psi_k(1)),\end{equation*}
where $\Psi_k(\psi_k(t))=y_k(t)$. So, we can conclude that $\Psi_k(x)= A_kJ_{ \frac{1}{2}}(x)+ B_k J_{ -\frac{1}{2}}(x)$ for $x\in (0,\psi_k(1))$ and $A_k,B_k \in \mathbb R$. Consider $J_{- \frac{1}{2}}(x)=\sqrt{\frac{2}{\pi x}}\cos x$, then $\chi(t)\sqrt{\frac{2}{\pi \psi_k(t)}}\cos \psi_k(t)= \sqrt{\frac{2}{\pi}\lambda_k\psi^\prime_k} \cos \psi_k(t)$ does not converge to $0$ as $t$ goes to $0$. Thus, by continuity at time $0$ and \eqref{bcode}, we reach the conclusion that $y_k(t) =A_kJ_{\frac{1}{2}}(\psi_k(t)),$ for all $ t\in [0,1]$.
\end{proof}
Now we use Cor.~\ref{cormainthm} to deduce the KL expansion of $Y_1$ and $Y_2$ by making use of the closed form expression $J_{\frac{1}{2}}(x)=\sqrt{\frac{2}{\pi x}}\sin x$. The boundary condition (\ref{bcode}) yields 
$$\chi(1)J_{\frac{1}{2}}\left(\frac{1}{\sqrt{\lambda_k}} \int_0^1 f(s)  \diff s\right)=0\quad\text{ or, } \sin \left(\frac{3[(|a|+b)^{\frac{1}{3}}-b^{\frac{1}{3}}]}{|a|\sqrt{\lambda_k}}\right) =0, \quad a\neq 0.$$
Therefore, for $Y_1(t)=B^{\tilde{c}}(\int_0^t(as+b)^{-\frac{4}{3}}\diff s)$, $a>0,$ the $k$-th largest eigenvalue is 
\begin{equation}\label{ev}
\lambda^{Y_1}_k = \frac{9}{a^2k^2\pi^2}\left[(a+b)^{\frac{1}{3}}-b^{\frac{1}{3}}\right]^2.
\end{equation}
The corresponding orthonormal eigenfunction is given by 
 \begin{equation}\label{Y1ef}[0,1]\ni t\mapsto \phi_k^{Y_1} (t)= \sqrt{\frac{2a}{3 [(a+b)^{\frac{1}{3}}-b^{\frac{1}{3}}]}} (at+b)^{-\frac{1}{3}} \sin \left(  \frac{k\pi [(at+b)^{\frac{1}{3}}-b^{\frac{1}{3}}] }{(a+b)^{\frac{1}{3}}-b^{\frac{1}{3}}}  \right).  \end{equation}
Thus, the KL expansion of $Y_1$ can be written as
\begin{align*}Y_1(t)= \sum_{k=1}^\infty  \xi_k \sqrt{\frac{6[(a+b)^{\frac{1}{3}}-b^{\frac{1}{3}}]}{ak^2\pi^2 }} (at+b)^{-\frac{1}{3}} \sin \left(  \frac{k\pi [(at+b)^{\frac{1}{3}}-b^{\frac{1}{3}}] }{(a+b)^{\frac{1}{3}}-b^{\frac{1}{3}}}  \right), \end{align*}
where $\xi_k \overset{\text{i.i.d.}}{\sim} \mathcal{N}(0,1)$.

For $ Y_2(t) =B^{\tilde{c}}(\int_0^t(|a|(1-s)+b)^{-\frac{4}{3}}\diff s$, $a<0$, we find that the eigenvalues have the expression given in \eqref{ev}, with $a$ replaced by $|a|$ and the orthonormal eigenfunction for each $k$ is 
$$t \mapsto \phi_k^{Y_2} (t) =\sqrt{\frac{2|a|}{3 [(|a|+b)^{\frac{1}{3}}-b^{\frac{1}{3}}]}}(|a|(1-t)+b)^{-\frac{1}{3}} \sin \left(  \frac{k\pi \left[ (|a|+b)^{\frac{1}{3}} - (|a|(1-t)+b)^{\frac{1}{3}} \right] }{[(|a|+b)^{\frac{1}{3}}-b^{\frac{1}{3}}]}\right).$$
The KL expansion of $Y_2$ is given by 
$$Y_2(t)= \sum_{n=1}^\infty  \xi_k \sqrt{\frac{6[(|a|+b)^{\frac{1}{3}}-b^{\frac{1}{3}}]}{|a|k^2\pi^2 }}  (|a|(1-t)+b)^{-\frac{1}{3}} \sin \left(  \frac{k\pi \left[ (|a|+b)^{\frac{1}{3}} - (|a|(1-t)+b)^{\frac{1}{3}} \right] }{(|a|+b)^{\frac{1}{3}}-b^{\frac{1}{3}}}  \right).$$
In the next section, we will utilize the formulas for eigenvalues and eigenfunctions derived above.
\section{Numerical Illustrations}\label{numericalsec}
\subsection{Eigenfunction and Variance}
\begin{figure}
\centering
\begin{subfigure}{0.49\textwidth}
\centering
    \includegraphics[width=\linewidth]{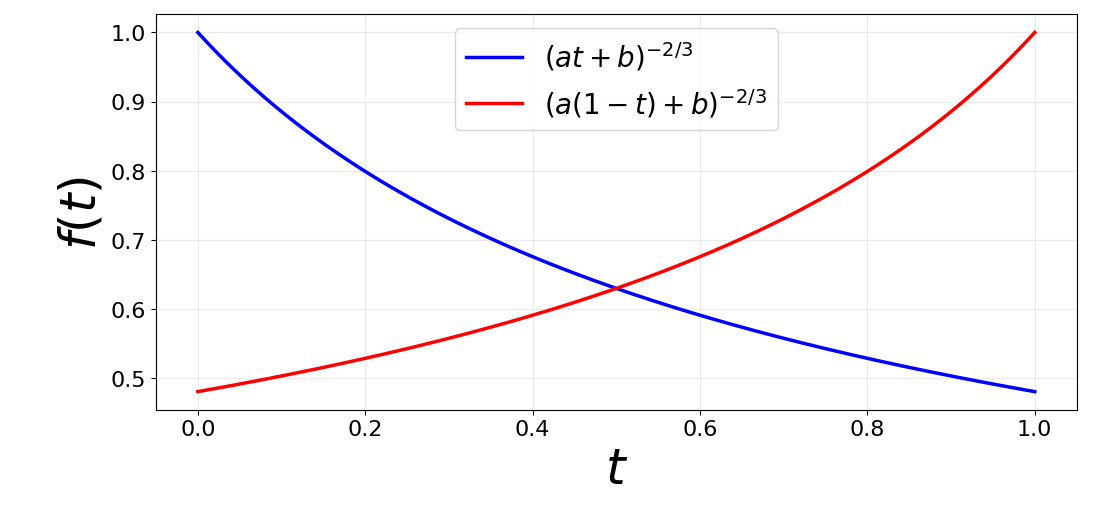}
    \caption{$f_1(t)$ (blue) and $f_2(t)$ (red) for $a=2, b=1$}  
    \label{f1f2}
\end{subfigure}
\begin{subfigure}{0.49\textwidth}
\centering
    \includegraphics[width=\linewidth]{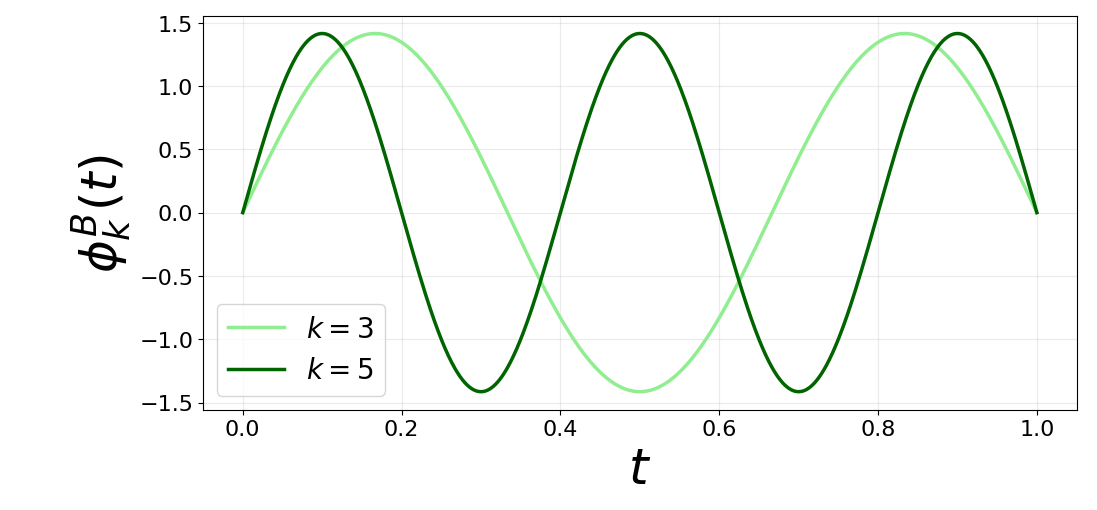}
    \caption{Eigenfunctions of $\mathscr{C}_B$}  
    \label{Bef}
\end{subfigure}
\begin{subfigure}{0.49\textwidth}
\centering
    \includegraphics[width=\linewidth]{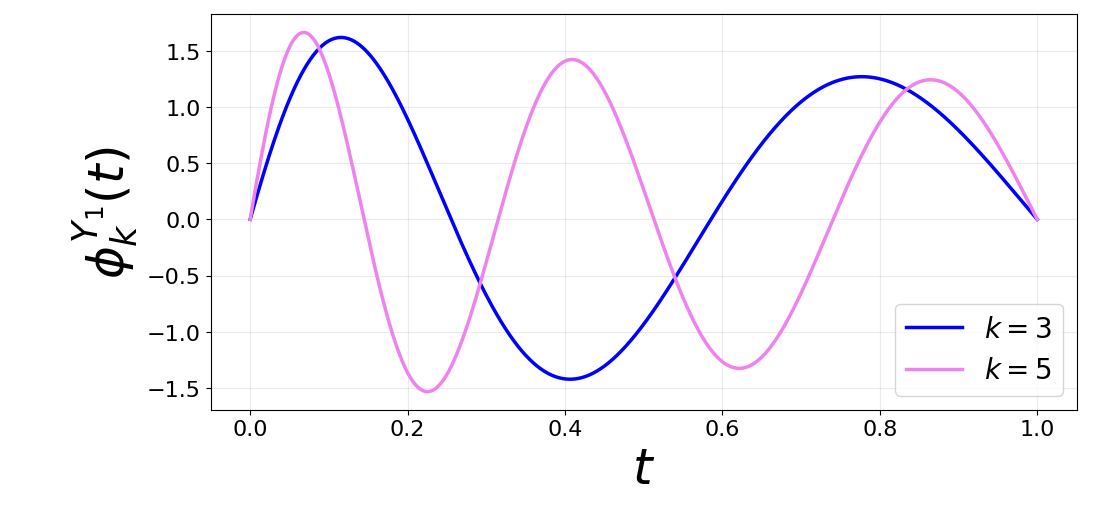}
    \caption{Eigenfunctions of $\mathscr{C}_{Y_1}$ for $a=2, b=1$} 
    \label{y1ef}
\end{subfigure}
\begin{subfigure}{0.49\textwidth}
\centering
    \includegraphics[width=\linewidth]{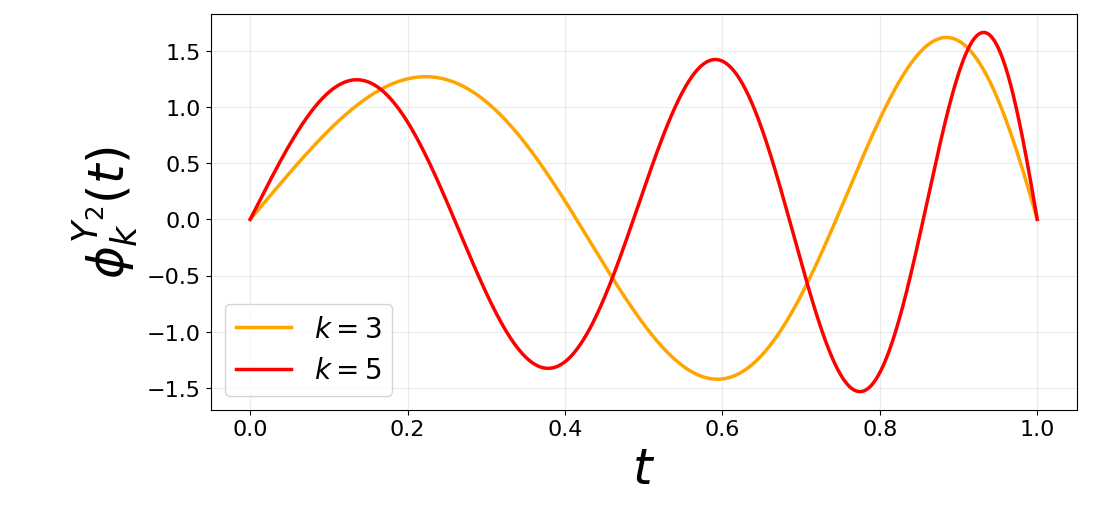}
    \caption{Eigenfunctions of $\mathscr{C}_{Y_2}$ for $a=-2, b=1$} 
    \label{y2ef}
\end{subfigure}
\caption{Illustration of asymmetry in time rescale speed and eigenfunctions of $\mathscr{C}_{Y_1}$ and $\mathscr{C}_{Y_2}$}
\label{figfunctions}
\end{figure}

\begin{figure}
\centering
\begin{subfigure}{0.49\textwidth}
\centering
    \includegraphics[width=\linewidth]{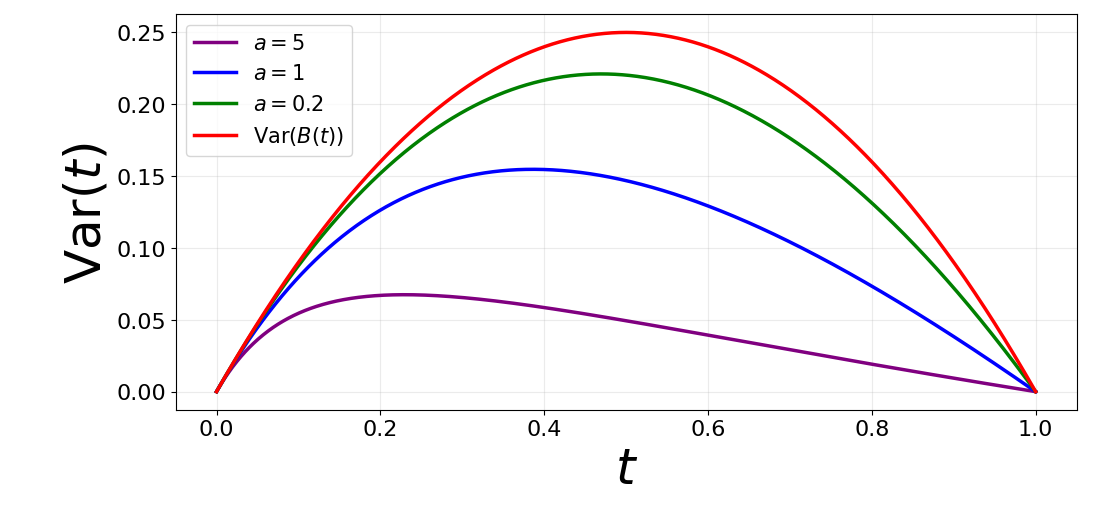}
    \caption{Variance function of $Y_1$ for fixed $b=1$}
    \label{Y1varwrta}
\end{subfigure}
\begin{subfigure}{0.49\textwidth}
\centering
    \includegraphics[width=\linewidth]{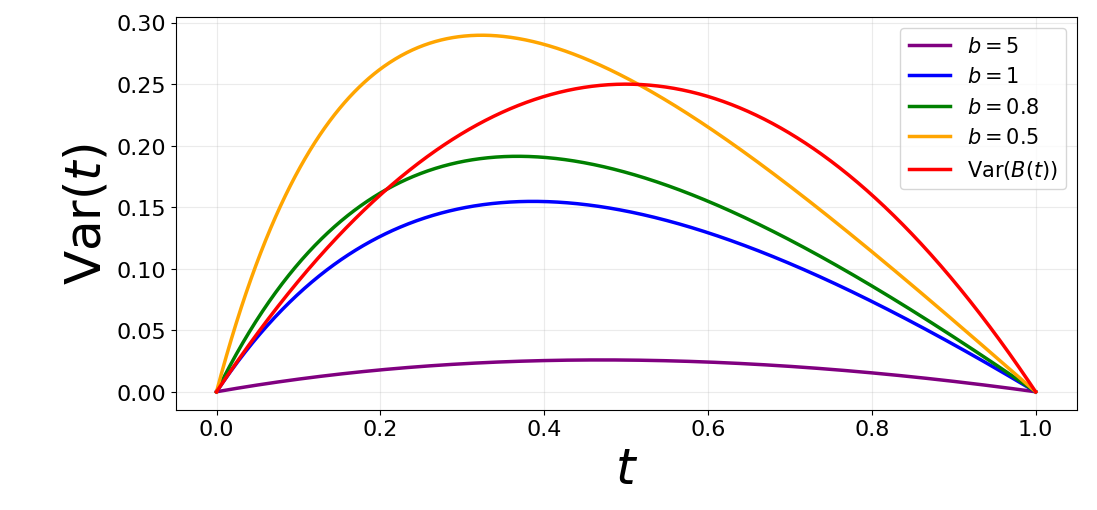}
    \caption{Variance function of $Y_1$ for fixed $a=1$} 
    \label{Y1varwrtb}
\end{subfigure}
\begin{subfigure}{0.49\textwidth}
\centering
    \includegraphics[width=\linewidth]{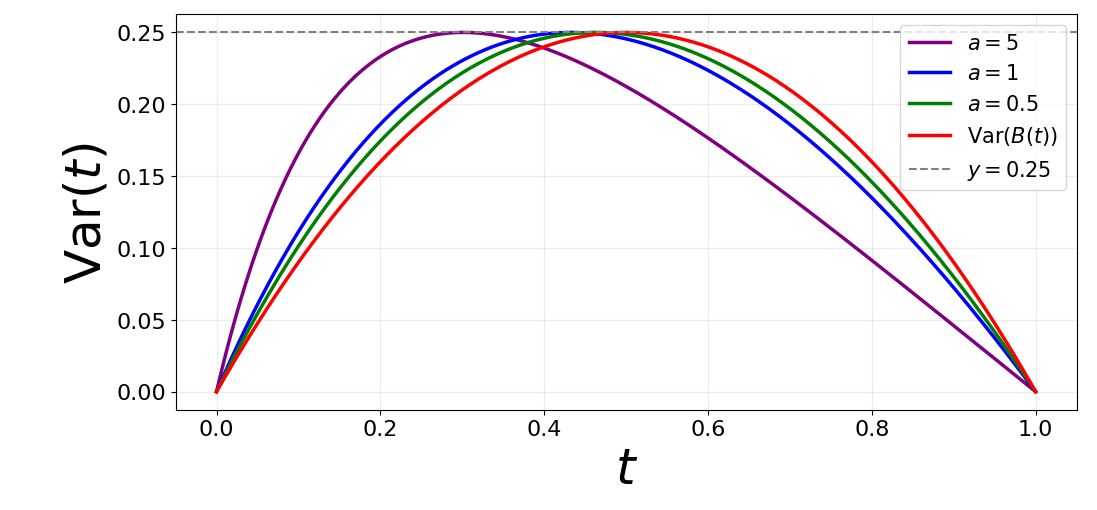}
    \caption{Variance function of $Y_{1,e}$ for fixed $b=2$} 
    \label{Y1evarwrta}
\end{subfigure}
\begin{subfigure}{0.49\textwidth}
\centering
    \includegraphics[width=\linewidth]{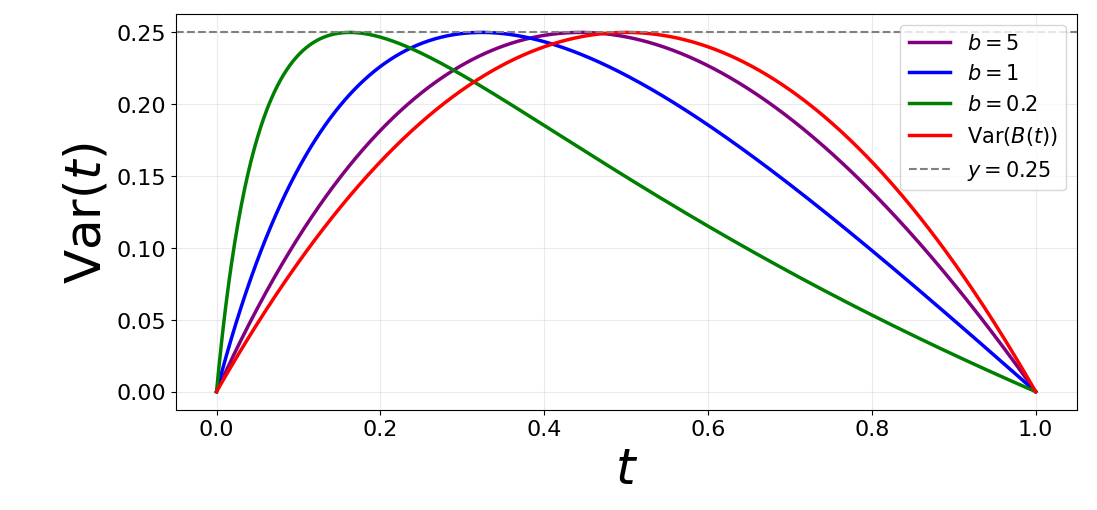}
    \caption{Variance function of $Y_{1,e}$ for fixed $a=2$}  
    \label{Y1evarwrtb}
\end{subfigure}
\caption{Variance analysis for $Y_1$}
\label{figY1}
\end{figure}

In this subsection, we illustrate the nature of asymmetry of the proposed two-parameter family of TRBB using their eigenfunctions and variance. Figs.~\ref{Bef}, \ref{y1ef}, and \ref{y2ef} show the few leading eigenfunctions of the standard Brownian bridge $B$ and the two TRBBs, $Y_1$ and $Y_2$. The eigenfunctions of $Y_1$ exhibit greater variability toward the left end of the interval, whereas those of $Y_2$ vary more toward the right end. The speed of time rescaling is measured by $\frac{\diff}{\diff t}\int_0^tf^2(s)\diff s=f^2(t)$, and its square root is plotted in Fig.~\ref{f1f2}. The numerical examples illustrate the monotonicity of the speed. Using \eqref{CVfnTRBB}, the variance function of $Y_1$, denoted as $V_1$, is given by 
$$V_{1}\left(t\right):=\operatorname{Var}(Y_1(t)) =\frac{3}{a}\left(b^{-\frac{1}{3}}-\left(at+b\right)^{-\frac{1}{3}}\right)\left(1-\frac{\left(b^{-\frac{1}{3}}-\left(at+b\right)^{-\frac{1}{3}}\right)}{\left(b^{-\frac{1}{3}}-\left(a+b\right)^{-\frac{1}{3}}\right)}\right),$$
where $a,b>0$. The speed of time rescaling is measured by $f_1^2(t)=(at+b)^{-\frac{4}{3}}$, which indicates that, for small values of $t$, time elapses much faster than for larger values of $t$; this explains the skewness of $V_1$ apparent in Figs.~\ref{Y1varwrta} and \ref{Y1varwrtb}. The TRBB $B^c(ct)$ on $[0,1]$ has a uniform time rescaling speed equal to $c$. Thus, if $a$ is small compared to $b$, the speed approaches uniformity, and the shape of $V_1$ approaches that of $B^{b^{-\frac{4}{3}}}(b^{-\frac{4}{3}}t)$, which is symmetric about $\frac{1}{2}$ and has maximum variance $\frac{b^{-\frac{4}{3}}}{4}$. This also explains the reduction in maximum variance and the increasingly symmetric shape as $b$ grows, as shown in Fig.~\ref{Y1varwrtb}.

\begin{figure}
\centering
\begin{subfigure}{0.49\textwidth}
\centering
    \includegraphics[width=\linewidth]{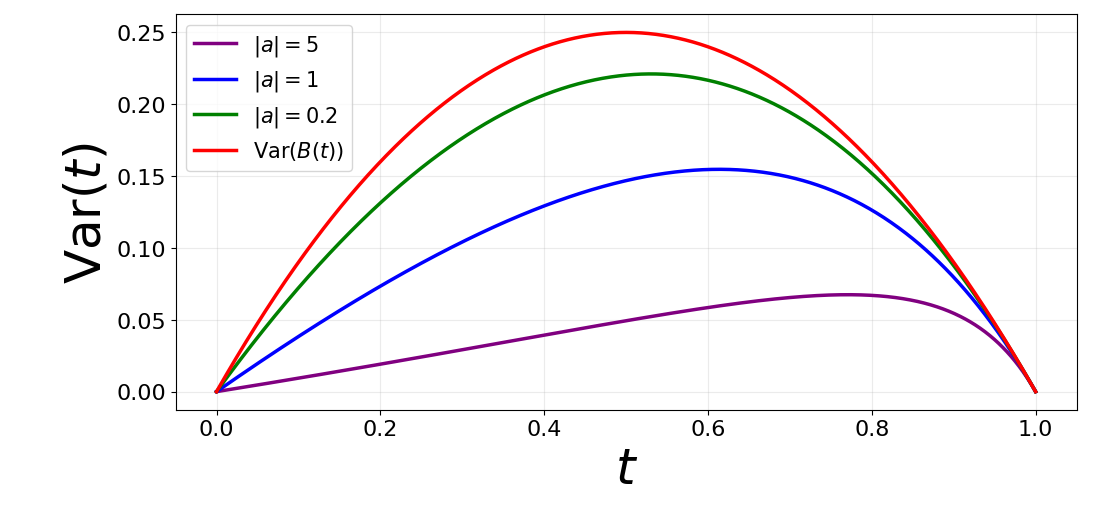}
    \caption{Variance function of $Y_2$ for fixed $b=1$}    
\end{subfigure}
\begin{subfigure}{0.49\textwidth}
\centering
    \includegraphics[width=\linewidth]{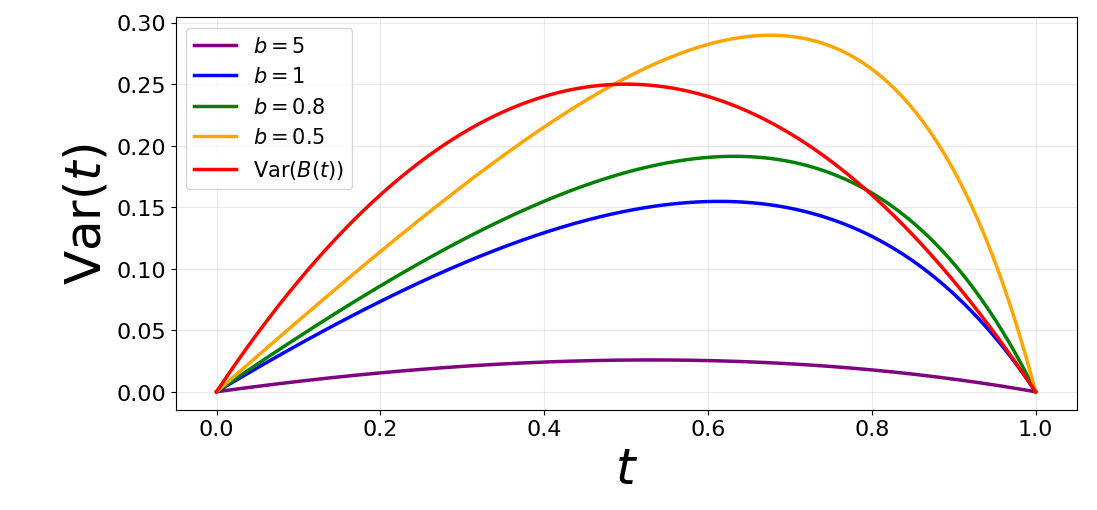}
    \caption{Variance function of $Y_2$ for fixed $a=-1$}    
\end{subfigure}
\begin{subfigure}{0.49\textwidth}
\centering
    \includegraphics[width=\linewidth]{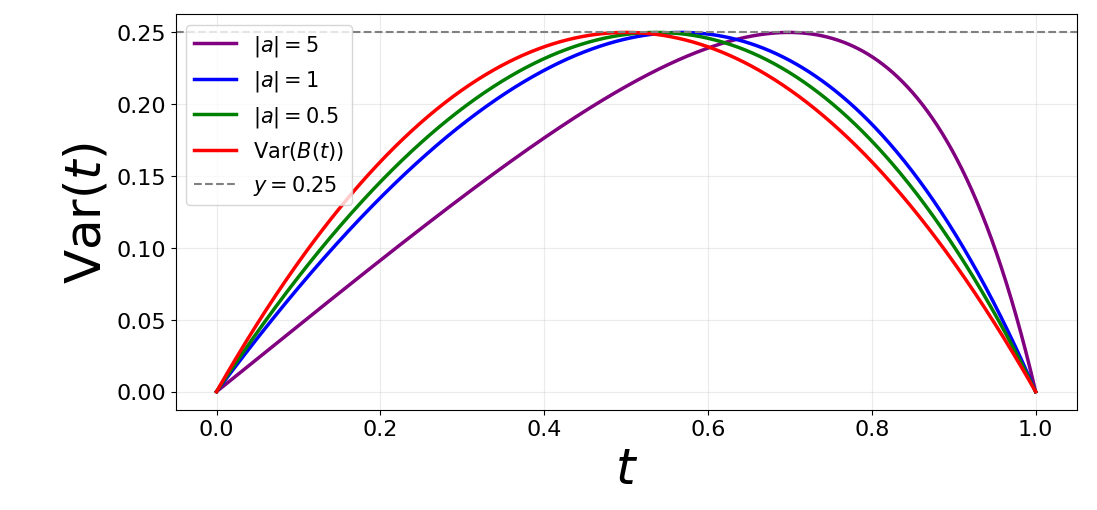}
    \caption{Variance function of $Y_{2,e}$ for fixed $b=2$}    
\end{subfigure}
\begin{subfigure}{0.49\textwidth}
\centering
    \includegraphics[width=\linewidth]{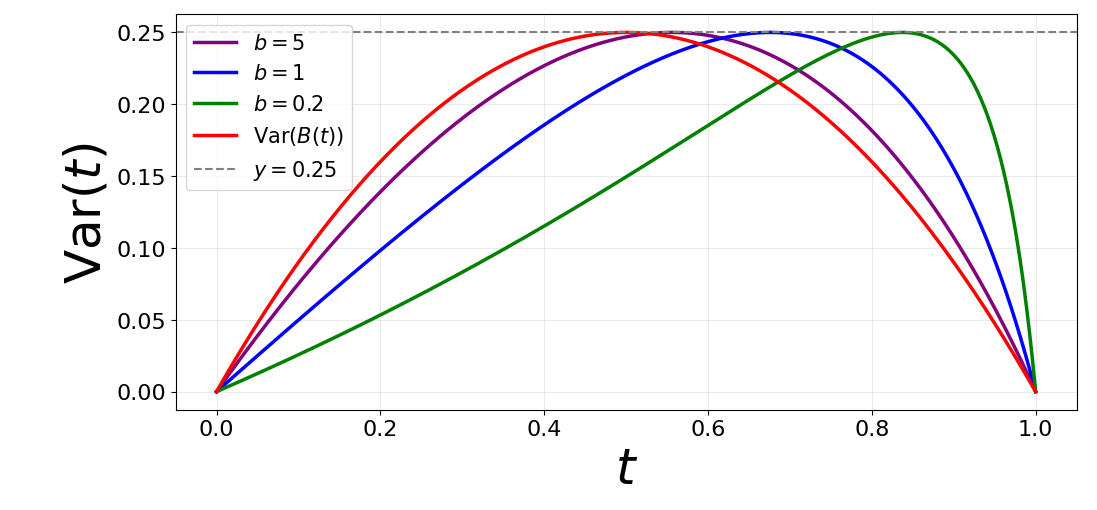}
    \caption{Variance function of $Y_{2,e}$ for fixed $a=-2$}    
\end{subfigure}
\caption{Variance analysis for $Y_2$}
\label{figY2}
\end{figure}

The three constants involved are $a$, $b$, and $\tilde c$. A formula for $b$ in terms of $a$ and $\tilde c$ is cumbersome. Therefore, to study the properties, we define $a_{e} :=\frac{a}{\tilde c^{-\frac{3}{4}}}$ and $b_{e} :=\frac{b}{\tilde c^{-\frac{3}{4}}}$. It turns out that $a_e$ and $b_e$, as functions of $a$ and $b$, respectively, are strictly increasing. Note that $c_e:=\int_0^1 (a_e t+b_e)^{-\frac{4}{3}}\diff t=1$. We define the corresponding TRBB as $Y_{1,e}(t)=B\left(\int_0^t(a_es+b_e)^{-\frac{4}{3}}\diff s\right)$, where $B$ is a standard Brownian Bridge on $[0,1]$. The variance function of $Y_{1,e}$ is shown in Figs.~\ref{Y1evarwrta} and ~\ref{Y1evarwrtb}. This shows how $a_e$ and $b_e$ affect the skewness of the variance function without changing its maximum value \eqref{CFBB}. We perform the same analysis as above for $Y_2$ as well. This is depicted in Fig.~\ref{figY2}.
\subsection{Empirical Power Comparison}\label{Power Comparison}
Here, the algorithm in Sec.~\ref{powersec} is used for the empirical power computation of the uniformity test against a specific non-uniform alternative with PDF $g^2$, where 
\begin{equation}\label{densityTRBB}g(t)=(at+b(a) - \min (0,a))^{-\frac{2}{3}}, \quad t\in [0,1].   \end{equation}
For each  $a\neq0$, $b(a)>0$ satisfies \begin{align}\label{b(a)}c(a,b):=\tilde{c}=\int_0^1(as+b(a)-\min(0,a))^{-\frac{4}{3}}\diff s =1,\end{align} 
so that $S$ is a CDF. Since $b\mapsto c(a,b)= \int_0^1 (as+b-\min(0,a))^{-\frac{4}{3}}\diff s$ is monotonic on $b\in (0,\infty)$, the implicit equation $c(a,b)=1$ admits a unique $b(a)$. For notational brevity, we write $b^*:=b(a)$. Using this notation, the CDF $S$ takes the form
\begin{equation}\label{TRBBCDF}S(t)= \int_0^t g^2(s)\diff s ={(b^*- \min (0,a))}^{-\frac{1}{3}}-(at+b^*- \min (0,a))^{-\frac{1}{3}}, \quad t\in [0,1].\end{equation}
The parameters for the experiment are $a=3.492$ and $b^*=0.168$, which are obtained from Table~\ref{tab:fit_summary_mle}.  Moreover, we consider the test statistic $s_n$ \eqref{TRBBstatistic} with $m=10$, and the significance level $\alpha=0.05$. Also, set $K_{\text{null}}=K_{\text{alt}}=20000$. The random seed is fixed (seed = 123456) and used to produce a pseudo-random number sequence using the code \verb|np.random.default_rng(seed)| in NumPy 1.26.4. 

\begin{figure}[t]
  \centering
  \includegraphics[width=0.85\textwidth]{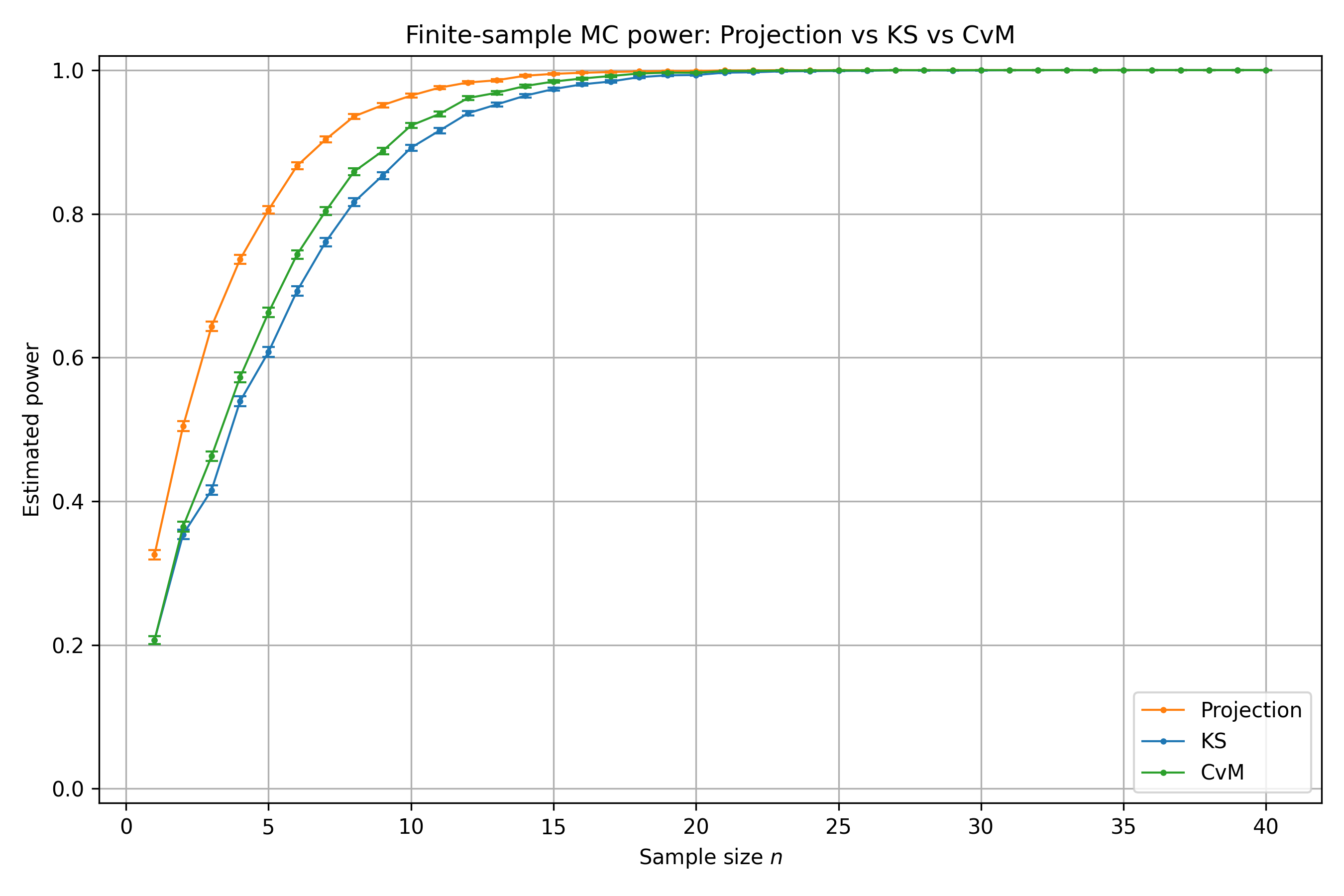}
  \caption[Estimated Power of KS, CvM, Projection Tests]{Estimated power of Projection, KS, and CvM tests using finite-sample MC. Vertical line segments denote approximate $95\%$ MC confidence intervals for the estimated power.}
  \label{power3tests_FSMC}
\end{figure}

The analytic alternative CDF\eqref{TRBBCDF} and its inverse are used to generate Monte Carlo (MC) samples. We compare the estimated powers of the KS, CvM, and projection tests. The finite-sample MC experiment (Fig.~\ref{power3tests_FSMC}) reveals three clear patterns. First, the projection test attains substantially higher power than the KS and CvM tests at small sample sizes; for instance, at $n=5$, the estimated powers are approximately $0.806$, $0.608$, and $0.663$ for the projection, KS, and CvM tests, respectively. Second, the CvM test consistently outperforms the KS test across the range of sample sizes; for example, at $n=3$, the estimated powers are approximately $0.463$ and $0.416$, respectively. Third, for $n\geq20$, all three tests have essentially reached power one. These findings have natural explanations: under the chosen alternative, the projection statistic aggregates signal in a small number of principal directions (the first $m$ basis functions), so even a small sample produces a nonnegligible mean shift of the statistic (the signal grows like $\sqrt{n}$), which the projection test is designed to detect efficiently. The KS statistic uses supremum deviation and is therefore less sensitive to smooth, distributed departures from uniformity; by contrast, the CvM statistic integrates squared deviations and so accumulates distributed signal more effectively than KS, explaining CvM's higher power here. Finally, with $K_{\text{alt}}=20000$, the MC standard error is at most $\approx 0.0035$, so the observed superiority at small $n$ is statistically significant. 
\section{Real-data analysis}\label{datasec}
In this section, we demonstrate the utility of CDF\eqref{TRBBCDF} as a model for real datasets with values in the closed interval $[0,1]$. To this end, we consider a sample of the proportion of French speakers in 2014 across 88 of nearly 195 countries. This dataset was previously presented and studied in \cite{frenchdataset} by fitting various models, including the Log-Shifted Gompertz (LSG) and Beta distribution. Here, we fit three models to the observed proportions data: the CDF given by \eqref{TRBBCDF}, the Beta distribution, and the LSG model \cite{frenchdataset}. We compute a variety of fit metrics, reported in Table~\ref{tab:fit_summary_mle}, including $L^2$ error of the CDF, classical goodness-of-fit tests (KS and CvM), the log-likelihood, and information criteria (AIC, BIC). All three models are fitted by \emph{maximum likelihood estimation} (MLE). For CDF\eqref{TRBBCDF}, the MLE procedure optimizes the free parameter $a$ and then computes $b^*$ via the model's normalization constraint \eqref{b(a)}. The Beta distribution was fitted with MLE on $[0,1]$ via \texttt{stats.beta.fit} with loc=0 and scale=1, yielding $(\alpha,\beta)$. Finally, the LSG parameters are obtained via a direct numerical MLE. 
\begin{table}[t]
\centering
\setlength{\tabcolsep}{6pt}
\small
\begin{tabular}{lccc}
\toprule
 & \textbf{CDF\eqref{TRBBCDF}} & \textbf{Beta($\alpha,\beta$)} & \textbf{LSG} \\
Parameters & $(\hat{a}=\num{3.492397},$ & $(\hat{\alpha}=\num{0.538151},$ & $(\hat{\alpha}=\num{1.934721},$ \\
 & $\hat{b}=\num{0.167802})$ & $\;\hat{\beta}=\num{1.857795})$ & $\;\hat{\beta}=\num{0.705283})$ \\ \hline
$L^2$ error & $\num{0.015843}$ & $\num{0.036}$ & {\bfseries \num{0.015066}} \\
KS statistic & {\bfseries \num{0.066680}} & $\num{0.096592}$ & $\num{0.073299}$ \\
KS $p$-value & {\bfseries \num{0.803965}} & $\num{0.361145}$ & $\num{0.703823}$ \\
CvM statistic & $\num{0.046577}$ & $\num{0.170591}$ & {\bfseries \num{0.043399}} \\
CvM $p$-value & $\num{0.897874}$ & $\num{0.332934}$ & {\bfseries \num{0.915982}} \\
log-likelihood & $\mathbf{58.5}$ & $54.1$ & $57.7$ \\
AIC & $\mathbf{-115}$ & $-104$ & $-111$ \\
BIC & $\mathbf{-112}$ & $-99$ & $-106$ \\
\bottomrule
\end{tabular}
\caption[Fit summary]{Fit summary (MLE). Columns report fitted parameter values, selected goodness-of-fit and information-criterion values for each model.}
\label{tab:fit_summary_mle}
\end{table}

\begin{figure}[t] 
  \centering
  \includegraphics[width=0.85\textwidth]{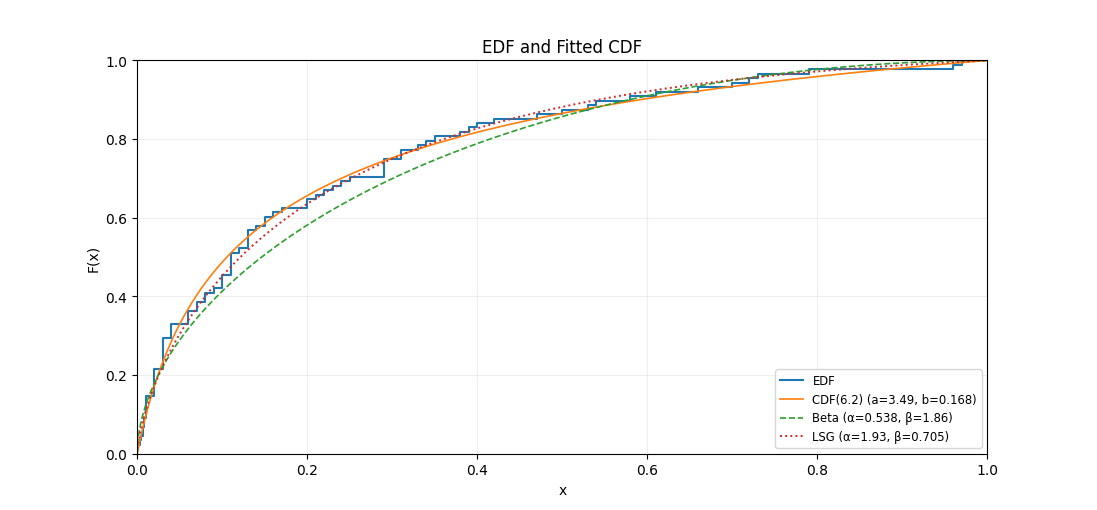}
  \caption[EDF and Fitted CDF]{EDF and Fitted CDF for the French speaking population data}
  \label{EDFandCDF}
\end{figure}

We calculate the $L^2$ error between the fitted CDF and the EDF using
$$\left(\int_0^1\big(F_{\text{model}}(t)- F_n(t)\big)^2\,\diff t\right)^{\frac{1}{2}},$$
where a smaller value indicates a better fit. The LSG model has the smallest $L^2$ error ($1.51\times10^{-2}$), while CDF\eqref{TRBBCDF} has a nearly identical value ($1.58\times10^{-2}$). Both models perform substantially better than the Beta distribution ($3.6\times10^{-2}$), indicating that the LSG model and CDF\eqref{TRBBCDF} track the EDF more closely than the Beta distribution (see Fig.~\ref{EDFandCDF}).

The computed KS statistics are modest for all three models. The CDF\eqref{TRBBCDF} (KS$=0.0667$, $p\approx0.804$) and the LSG model (KS$=0.0733$, $p\approx0.704$) yield large $p$-values, so we fail to reject the goodness-of-fit null hypothesis at conventional significance levels. Although the Beta distribution (KS $=0.0966$, $p\approx0.361$) has a larger KS statistic, it is still not rejected at conventional significance levels. 

The CvM statistic is smallest for the LSG model, followed closely by CDF\eqref{TRBBCDF}, and substantially larger for the Beta model. This ordering indicates that LSG and CDF\eqref{TRBBCDF} match the EDF more closely than the Beta model. The corresponding CvM $p$-values are such that none of the models is rejected at common significance levels.

CDF\eqref{TRBBCDF} attains the largest log-likelihood ($58.5$), followed closely by the LSG model ($57.7$), while the Beta model has a substantially lower log-likelihood ($54.1$). CDF\eqref{TRBBCDF} also has the smallest (best) values for the penalized criteria (AIC $=-115$, BIC $=-112$), followed by the LSG model.
\begin{remark}
Under the implemented MLE procedure, CDF\eqref{TRBBCDF} and the LSG model provide markedly better fits to the EDF than the Beta model, with the LSG model performing slightly better in terms of raw CDF-distance ($L^2$ error and CvM), while CDF\eqref{TRBBCDF} achieves the highest log-likelihood and the best information-criterion values. The Beta model performs less well across multiple metrics. Comparing our results for CDF\eqref{TRBBCDF} with the corresponding results for other competing models reported in \cite{frenchdataset,3datasets} shows that the proposed CDF achieves a better fit than these models.
\end{remark}
\begin{remark}
Using the fitted CDF\eqref{TRBBCDF} for parameter estimates $\hat{a}=\num{3.492397}$ and $\hat{b}=0.168$, we compute the projection statistic $s_n$ \eqref{TRBBstatistic} with $m=10$. For the given dataset, we obtain $s_n=0.145$. Under the hypothesis that the population follows the fitted CDF\eqref{TRBBCDF}, the statistic $s_n$ is asymptotically standard normal; hence, the corresponding two-sided asymptotic $p$-value is $88.5\%$. Since this $p$-value is large, the fitted model explains the data satisfactorily. 
\end{remark}
\section{Conclusion}\label{conlsec}
We propose a goodness-of-fit statistic \eqref{TRBBstatistic} for unit distributions, constructed by projecting the empirical process onto the principal eigenfunctions of the covariance operator of TRBB \eqref{UDTRBB}. Moreover, for the uniformity test \eqref{H_0}, using \eqref{TRBBstatistic}, against a known density as the alternative, an asymptotic distribution of its power is derived. We illustrate that the power for the specific case \eqref{densityTRBB} is higher than KS and CvM tests, for the parameter values considered.

Furthermore, we derived the KL expansion of the Brownian bridge $B^{\tilde c}$ with the time rescaling $t\mapsto \int_0^t f^2(s)\diff s$ for the two-parameter family $f(s)= (as +b - \min (0,a)\;)^{-\frac{2}{3}}$, $a\neq 0$, $b>0$, where $\tilde c=\int_0^1 f^2(s)\diff s$ and $B^{\tilde c}$ is the standard Brownian bridge on $[0,{\tilde c}]$ \eqref{BB}. Due to the non-uniformity of $f$, the variance function of the proposed TRBB family is skewed, unlike that of the standard Brownian bridge. The eigenfunctions of the covariance operator are tractable and are expressed using the sine function.  Finally, we show that the associated CDF\eqref{TRBBCDF} is a competitive model for a real dataset and has better goodness-of-fit measure values when compared with two other models.  

As per our knowledge, there are not many results for multiparameter Brownian bridges and their time rescaling in the literature. A well-chosen diverse class of multiparametric rescaling yielding a tractable KL basis is worth seeking, which could be a future research direction. Apart from its theoretical interest, such investigations are also interesting for their potential application in the inference problem for a multivariate distribution with compact support. 
\section{Acknowledgements}
The second author acknowledges the financial support from the National Board for Higher Mathematics (NBHM), Department of Atomic Energy, Government of India, under grant number 02011-32-2023-R\&D-II-13347. The third author acknowledges the financial support received through the Kishore Vaigyanik Protsahan Yojana Fellowship (20092903/X063020220) during 2021–25. The third author also acknowledges the financial support of the Student Research Award Programme by the Centre for Intelligent Solutions, IISER Pune (2025–26), which supported the research collaboration with G2Q Computing SRL, Italy.
\section{Declaration}
\begin{itemize}
\item Conflict of interest: The authors declare no conflict of interest.
\item AI tools: The draft code for the numerical experiment was generated by ChatGPT 5.3 and was thoroughly reviewed and rectified before execution. The source code can be found in the GitHub repository \cite{Nikhil2026github_TRBB}.
\item Availability of data and materials: Data sharing is not applicable to this article as no new data was created or analyzed in this study. All empirical studies reported in this article are based on openly available data.
\item Authors’ contributions: Authors’ names appear in alphabetical order of their surnames. The first author participated in formulating applications of theoretical results, improving presentations, and securing funding. The second author participated in revising research questions, cross-checking all theoretical results, and drafting and revising of the manuscript. The third author was responsible for ideation, formulation of the methodology, theoretical development, and drafting and revising of the manuscript.
\end{itemize}

\nocite{NikhilMSthesis}
\printbibliography
\end{document}